\documentclass[a4paper]{article} 
\usepackage{graphicx} 

\long\def\comment#1{}

\usepackage{a4wide,amsmath}
\makeatletter
\def\@normalsize{\@setsize\normalsize{10pt}\xpt\@xpt
\abovedisplayskip 10pt plus2pt minus5pt\belowdisplayskip
\abovedisplayskip \abovedisplayshortskip \z@
plus3pt\belowdisplayshortskip 6pt plus3pt
minus3pt\let\@listi\@listI}

\def\subsize{\@setsize\subsize{12pt}\xipt\@xipt}
\def\section{\@startsection {section}{1}{\z@}{1.0ex plus
1ex minus .2ex}{.2ex plus .2ex}{\large\bf}}
\def\subsection{\@startsection
   {subsection}{2}{\z@}{.2ex plus 1ex} {.2ex plus .2ex}{\subsize\bf}}
\makeatother

\begin{document}

\date{}

\title{\huge \bf {Equivalence of the categories of group triples and of hypergroups over the group}}


\vspace{10mm}

\author{\bf{Samuel Dalalyan}\\\\\vspace{5pt}
Department of Mathematics and Mechanics,\\
Yerevan State University, Armenia 
}

\maketitle

\begin{abstract}
The main result of this paper is that the categories of (right) hypergroups over  the group and of triples,
 consisting of a group, its subgroup and a (right) transversal to this subgroup,
 are equivalent.\\\\\vspace{5pt}
\end{abstract}
  

{\bf Introduction.}
 

\hspace{10mm} 
The concept of a (right) hypergroup over  the group was introduced in \cite{D1} and was refined in \cite{D2} and \cite{D3}. 
This concept generalizes and  unifies the concepts of the group, of the  field and of the linear space over the  field. 
The hypergroups over the  group are  already successfully applied to the generalization 
of the Schreier theorem about group extensions (\cite{D3}), in the theory of group  cohomology 
and in some other areas.

\hspace{10mm} 
All (right) hypergroups over the group and their morphisms determine a category, which is denoted ${\mathcal H}g$.
On the other hand, we consider the  category ${\mathcal G}Trip$ of (right) group triples and their morphisms. 
The objects of this category are  (right) group triples $(G, H, M)$, where $G$ is a group, 
$H$ is a subgroup of $G$ and $M$ is a (right) transversal to the subgroup $H$. 
The group triples generalizes the short exact sequence of groups 
together with a fixed section. 
The main result of this article is the following theorem.

\hspace{10mm}
{\bf Theorem 1.}
{\it The categories ${\mathcal H}g$ and ${\mathcal G}Trip$ are equivalent.}

\hspace{10mm} 
To prove this theorem, we need to construct two functors 
\begin{center}
${\bf H}:  {\mathcal G}Trip \rightarrow  { \mathcal H}g$ and 
${\bf T}: {\mathcal H}g \rightarrow {\mathcal G}Trip $
\end{center}
such that 
the composites ${\bf H} \circ {\bf T}$ and ${\bf T} \circ {\bf H}$ 
are naturally isomorphic to the identity functors 
$1_{{\mathcal G}Trip}$ 
and $1_{{\mathcal H}g}$.

\hspace{10mm} 
The present paper consists of this Introduction and five sections. 
In section 1 the category of (right) hypergroups over the group 
and in  section 2 the category of (right) group triples are defined.
In  section 3 the functor $\bf H$ and 
in  section 4 the functor $\bf T$ are constructed. 
In section 5 the main results of this paper (Theorem 1) is proved.


\hspace{10mm} 
{\bf Remark.}
Most recently, after the writing of this article, the author became aware of the article \cite{L}.
The paper \cite{L} is written on the basis of the same ideas, as our paper.
However it must be emphasized, that the main object entered into circulation in \cite{L} 
and called   the ${\it c-groupoid}$, is  a special case of the hypergroup over the group.
The notion  c-group coincides with the our notion  {\it unitary} hypergroup over the group (\cite{DN}).
Accordingly, all  results of \cite{L} follow from  the corresponding results of the papers \cite{D2}, \cite{D3}. 
In addition, I would like to especially note that we use  a notation system, 
which helps to better represent and memorize all relationships in this new theory.

	
\vspace{3mm}

{\bf 1.  The category of hypergroups over the group.}
	Let $M$ be a set,  the elements of  $M$ will be denoted by small Latin letters, 
	$H$ be a group, the elements of $H$ will be denoted by small Greek letters. 
	Consider the Cartesian products $M \times H$ and $M \times M$, 
	and four mappings
	$$
\Phi: M \times H \rightarrow M,  \qquad
\Psi: M \times H \rightarrow H,  \qquad
\Xi: M \times M \rightarrow M,  \qquad
\Lambda: M \times M \rightarrow H,  
$$
for which the following denotations are used:
$$
	\Phi(a,\alpha) = a^\alpha, \quad
	\Psi(a,\alpha) = {^a}\alpha, \quad
	\Xi(a, b) = [a, b], \quad
	\Lambda(a, b) = (a, b),
	$$
	and let  $\Omega = (\Phi, \Psi, \Xi, \Lambda)$.

	\hspace{10mm}
Consider the following conditions on mappings of this system.

\hspace{10mm} 
P1) The mapping $\Xi$ is a binary operation on $M$ with
the  properties
	
	(i)  any equation $[x, a] = b$ with elements $a, b\in M$ 
	has a unique solution in $M$;
	
	(ii) there exists a left neutral element $o \in M$, i.e. $o$ satisfies the condition  
	$[o, a] = a$ for any element $a \in M$.

	\hspace{10mm} 
	P2) The mapping $\Phi$ is a (right) action of  the group $H$ on $M$, that is
	
	(i) $(a^{\alpha})^{\beta}= a^{\alpha \cdot \beta}$ for any elements $\alpha, \beta \in H$ and $a \in M$;
	
	(ii) $a^{\varepsilon}= a$ for any $a \in M$, where $\varepsilon$ is the neutral element of the group $H$.
		
\hspace{10mm} 
P3) The mapping $\Psi$ sends the subset $\{ o \} \times H$ of $M \times H$ on $H$.  

\hspace{10mm} 
P4) The mappings of the system $\Omega$ satisfy the following identities:
	\begin{itemize}
	\item{$(A1)$} \quad  ${^a}(\alpha \cdot \beta) = {^a}\alpha \cdot {^{a^\alpha}}\beta$,
	\item{$(A2)$} \quad  ${[a, b]}^\alpha = [a^{{^b}\alpha}, b^\alpha]$,
	\item{$(A3)$} \quad  $(a, b) \cdot {^{[a,b]}}\alpha = {^a}({^b}\alpha) \cdot (a^{{^b}\alpha}, b^\alpha)$,
	\item{$(A4)$} \quad  $[[a, b], c] = [a^{(b,c)}, [b, c]]$,
	\item{$(A5)$} \quad  $(a, b) \cdot ([a, b], c) = {^a}(b,c) \cdot (a^{(b, c)}, [b, c])$.
	\end{itemize}

\hspace{10mm} 
{\bf Definition 1.1.}	A (right) {\it hypergroup over the group} is an object, which is defined by a database, consisting of

- a ({\it basic}) set $M$, 

- a group $H$, 

- a system of ({\it structural}) mappings  $\Omega = (\Phi, \Psi, \Xi, \Lambda)$, 

for which the conditions P1) - P4) are  satisfied. 
Such a  hypergroup over the group is denoted by $M_H$. 
Dually a left hypergroup over the group  is defined.

\hspace{10mm} 
	In algebra the term hypergroup is already used in an entirely different sense 	(see \cite{Marty}, \cite{Wall}, also \cite{Lit}). 
	In this paper only the right hypergroups  $M_H$ over the group are used 	and they are often shortly named  hypergroups.
	\hspace{10mm}


\hspace{10mm}
{\bf Proposition 1.1.} {\it Any hypergroup $M_H$ has the following additional properties:
\begin{itemize}
\item{$(A6)$} $\quad ^a\varepsilon = \varepsilon$;
\item{$(A7)$} $\quad o^\alpha = o$;
\item{$(A8)$} $\quad ^o\alpha = \theta^{-1} \cdot \alpha \cdot \theta$;
\item{$(A9)$} $\quad  (o, a) = \theta^{-1}$,
\item{$(A10)$} $\quad [a, o] = a^{\theta^{-1}}$;
\item{$(A11)$} $\quad (a, o) = {^a}(\theta^{-1})$.
\end{itemize}
 where $a \in M, \alpha \in H$ are arbitrary elements, $o$ is the left neutral element of the binary operation $\Xi$,
$\varepsilon$ is the neutral element of the group $H$ and  $\theta = \Lambda(o, o)^{-1}$.}


\hspace{10mm} 
{\bf Definition 1.2.}
Let $M_H$ and $M'_{H'}$ be hypergroups with systems of structural mappings 
$\Omega = (\Phi,\Psi,\Xi,\Lambda)$ and $\Omega' = (\Phi',\Psi',\Xi',\Lambda')$, respectively. 
A {\it morphism} 
$$
f: M_H \rightarrow M'_{H'}
$$ 
of hypergroups over the  group is a pair $f = (f_0, f_1)$, 
consisting of a homomorphism of groups 
$f_0: H \rightarrow H'$ 
and of a map of sets 
$f_1:M \rightarrow M'$, 
preserving the structural mappings, i.e. satisfying the following relations:
\begin{itemize}
\item{$(M\Phi)$} \quad  $\Phi \circ f_1 = (f_1\times f_0)\circ \Phi'$, 
\item{$(M\Psi)$} \quad  $\Psi \circ f_0 = (f_1\times f_0) \circ \Psi'$, 
\item{$(M\Xi)$} \quad  $\Xi \circ f_1 = (f_1\times f_1) \circ \Xi'$, 
\item{$(M\Lambda)$} \quad  $\Lambda \circ f_0 = (f_1\times f_1) \circ \Lambda'$.
\end{itemize}

\hspace{10mm} 
{\bf Remark 1.1.}
The condition $(M\Xi)$ means that $f_1$ is a homomorphism from $(M, \Xi)$ to $(M', \Xi')$. 
It is not difficult  to check that if  $(M, \Xi)$ and $(M', \Xi')$ satisfy the condition P1) with left neutral elements $o, o'$, respectively, 
and $f_1$ is a homomorphism from $(M, \Xi)$ to $(M', \Xi')$, then $f_1(o) = o'$. 

\hspace{10mm} 
{\bf Remark 1.2.} 
Any morphism $f = (f_0, f_1): M_H \rightarrow M'_{H'}$ sends the element $\theta = \Lambda(o, o)^{-1}$ 
to the element $\theta' = \Lambda(o', o')^{-1}$, where $o$ and $o'$ are left neutral elements of $(M, \Xi)$ 
and $(M', \Xi')$, respective,y. This immediately follows from $(M\Lambda)$.

\hspace{10mm} 
The simplest example of a morphism of hypergroups is the {\it identity morphism} of an hypergroup $M_H$:
$$
1_{M_H} = (1_H, 1_M): M_H \rightarrow M_H,
$$
 where $1_H$ and $1_M$ are the identity maps of $H$ and $M$. 
For arbitrary morphisms of hypergroups
$$
f = (f_0, f_1) :M_H \rightarrow M'_{H'}, \qquad f' = (f'_0, f'_1): M'_{H'} \rightarrow M''_{H''}
$$
the {\it composite} $f  \circ f' : = (f_0 \circ f'_0, f_1 \circ f'_1)$ is a morphism of hypergroups, as well.
The composite of morphisms of hypergroups satisfies the associative law:
$$
(f \circ f') \circ f'' = f \circ (f' \circ f''),
$$ 
and the neutrality law for identical morphisms is true:
$$
1_{M_H} \circ f = f = f \circ 1_{M'_{H'}}.
$$
 Thus, the classes of hypergroups over the group and their morphisms determine a category,  
 {\it the category of (right) hypergroups over the group}, which is denoted by ${\mathcal H}g$.


	\hspace{10mm} 
	{\bf Example 1.1.} Any group $M$ can be considered as a  hypergroup over the trivial group $H = \{ \varepsilon \}$. 
	Then
	
	- there exist, evidently,  unique (and {\it trivial}) mappings $\Psi$ and $\Lambda$, 
	
	- the mapping $\Phi$, satisfying P2(ii), also determined uniquely. 
	
Concerning $\Xi$,	we take it coincide with the binary operation of the group $M$. 
	This system of mappings $\Omega$ will satisfy the conditions P1) - P4). 
	Consequently, we obtain a hypergroup $M_{ \varepsilon \}}$. 
	
	For any homomorphism of groups $f_1: M \rightarrow M'$, 
	a (unique) morphism of hypergroups $f = (f_0, f_1):  M_{\{ \varepsilon \}} \rightarrow M'_{\{  \varepsilon \}}$ is determined.
	
	Thus, we obtain a natural embedding of the category $\mathcal G$ of groups  into the category ${\mathcal H}g$.
	
	\hspace{10mm}
	Emphasize that for any hypergroup $M_H$ over the trivial group  $H = \{ \varepsilon \}$, 
	we have that $(M, \Xi)$ is a group.
	Indeed, in this case according to (A4) and P2) (ii) the binary operation $\Xi$ is associative,
	and this assertion follows from the well known facy, that
	any binary operation, which satisfies the conditions  P1) and the  associative law, is a group operation.

		\hspace{10mm} 
	{\bf Example 1.2.} A hypergroup over the group can be canonically associated to any field $k$ in a following way.
	Let $M$ be the additive group, $H$ be the multiplicative group of $k$.
	Define the system of structural mappings $\Omega = (\Phi, \Psi, \Xi, \Lambda)$ in such a manner: 
	
	- $\Xi$ is the addition operation of $M$; 
	
	- $\Phi(a, \alpha)$ is  the product $a \cdot \alpha$ in $k$; 
	
	- $\Psi(a, \alpha) = \alpha$ for any $a \in M, \alpha \in H$;
	
	- $\Lambda(a, b) = \varepsilon$, where $a, b \in M$ and $\varepsilon$ is the neutral element of $H$. 
	
	These mappings satisfy the conditions P1) - P4), consequently, we get a hypergroup over the group   $M_H$, 
	{\it associated with the field $k$}.
	
	Let $f_1: k \rightarrow k'$ be a (mono)morphism of fields, 
	$f_1: M \rightarrow M'$ be the corresponding monomorphism of additive groups, 
	$f_0: H \rightarrow H'$ be the corresponding monomorphism of  multiplicative groups of these fields. 
	Then  $f = (f_0, f_1): M_H \rightarrow M'_{H'}$ gives a morphism of hypergroups over the group.
		Therefore there is  a natural embedding from the category $\mathcal F$ of fields 
		into the category ${\mathcal H}g$ of  hypergroups over the group.

		\hspace{10mm} 
	{\bf Example 1.3.}
	Similarly by considering for any linear space $L$ over the field $k$ 
	a hypergroup $M_H$, where 
	$M$  is the additive group of $L$, 
		$H$ is the multiplicatove group of the field $k$,
		a maturel embedding from the category $\mathcal L$ of lunear spaces over the field 
		into the category ${\mathcal H}g$ is determined.


	\vspace{3mm}
	
	{\bf 2. The category of group triples.}
Let $G$ be a group, $H$ be its subgroup.
A subset $M$ of $G$ is called 

(rt) a {\it right transversal} to the subgroup $H$ if	
any coset $Ha = \{ x \cdot a, x \in H \}$, $a \in G$ has a unique common element with the set $M$; 

(rcs) a {\it right complementary set}	to the subgroup $H$ if
for any element $x \in G$ there exists a unique representation 
$x = \alpha \cdot a$, where  $\alpha \in H$, $a \in M$.
	
	This two conditions on the subset $M \subset G$ are equivalent (see, for example,  \cite{R}).

\hspace{10mm} 
{\bf Definition 2.1.} 
{\it A right group triple} $(G, H, M)$ is a triple, consisting of a group $G$, of a subgroup $H$ of  $G$ and of a 
right transversal $M$ to the subgroup $H$.
Dually a left group triple is defined.

\hspace{10mm} 
A group triple can be considered as is a generalization of a short exact sequence of groups
$$
E \longrightarrow H \longrightarrow G \longrightarrow M \longrightarrow E, 
$$
together with a fixed section $\sigma: M \rightarrow G$. 
Such an object  we denote by $(\varphi, \psi, \sigma)$. 
The image $H'$ of $\varphi$ is a normal subgroup of $G$,
the image $M'$ of $\sigma$ is simultaneously a left and a right transversal to $H'$. 
Thus, we get  a  group triple $(G, H', M')$, which is simultaneously left and right group triple.

	\hspace{10mm} 
	Further we consider only the right group triples 
	and omit the word right.

	\hspace{10mm} 
	{\bf Definition 2.2.} A {\it morphism} of group triples
$$
g: (G, H, M) \rightarrow (G', H', M')
$$ 
 is a group homomorphism $g: G \rightarrow G'$ such that 
$$
g(H) \subset H', \quad  g(M) \subset M'.
$$ 
The composite of morphisms $g$ and $g': (G', H', M') \rightarrow (G'', H'', M'')$ is defined by
the composite $g \circ g'$ of corresponding group homomorphisms. 
The identity morphism $1_{(G, H, M)}$ is determined by the identity map $1_G$. 
The  class of all group triples together with all their morphisms determines a category,
the {\it category of group triples} which is denoted ${\mathcal G}Trip$.

Note that the class of all short exact sequences of groups with a section together with all their morphisms 
determines  a full subcategory of the category ${\mathcal G}Trip$.


	\vspace{3mm}
	{\bf 3. A standard construction for obtaining hypergroups over the group.}
	There is a standard method for obtaining hypergroups over the group.
	Let $(G, H, M)$  be an arbitrary group triple. 
	Since $M$ is a right complementary set to the subgroup $H$ of  $G$,  
	for any elements $a,b \in M$ and $\alpha \in H$ the elements $a \cdot \alpha$ and $a \cdot b$ 
	are uniquely represented as products of elements of $H$ and $M$. 
	Consequently, one can define  
	$$
	\Phi(a,\alpha) = a^\alpha, \quad \Psi(a,\alpha) = {^a}\alpha, \quad \Xi(a,b) = [a, b], \quad  \Lambda(a,b) = (a, b)
	$$  
	by using the relationships (St1) and (St2):
	\begin{itemize}
	\item{$(St1)$} \quad  $a \cdot \alpha = {^a\alpha} \cdot a^\alpha, \quad
	^a\alpha \in H, \,  a^\alpha \in M$,
	\item{$(St2)$} \quad  $a \cdot  b = (a, b) \cdot [a,b], \quad (a,b) \in H, \, [a,b] \in M$.
	\end{itemize}
	Thus we have a system of mappings $\Omega = (\Phi, \Psi, \Xi, \Lambda)$, canonically associated 
	to any group triple $(G, H, M)$.

		\hspace{10mm}	
	{\bf Proposition 3.1.} {\it For any group triple $(G, H, M)$ the canonically associated 
	system of mappings $\Omega = (\Phi, \Psi, \Xi, \Lambda)$ satisfies the properties P1)-P4), 
	consequently, determines a (right) hypergroup over the group $M_H$.}

It is said that this hypergroup $M_H$ over  the group {\it is obtained from the group triple $(G, H, M)$ 
by the standard construction}.	
	
		\hspace{10mm}	
{\bf Proof.}
	According to associative law of group binary operation we have
$$
(a \cdot \alpha) \cdot \beta = a \cdot (\alpha \cdot \beta), \quad 
(a \cdot b) \cdot \alpha = a \cdot (b \cdot \alpha), \quad
(a \cdot b) \cdot c = a \cdot (b \cdot c)
$$
for every $a, b, c \in M, \, \alpha, \beta \in H$. 
Applying the relations (St1), (St2) and the property (cs) of  uniqueness for the representation 
any element of $G$ as a product of elements from $H$ and $M$, we get, respectively, the pairs of relations 
P2(i) and (A1),  (A2) and (A3), (A4) and (A5).  

	\hspace{10mm}	
Similarly, using the equalitiy $a \cdot \varepsilon = \varepsilon \cdot a$, where $a \in M$ and $\varepsilon$ is 
the neutral element of $H$ (and consequently, of $G$),
we get the relations P2(ii) and (A6).

	\hspace{10mm}	
	Let the elements $\theta \in H, o \in M$ are uniquely determined by the relation $\varepsilon = \theta \cdot o$. 
	Then using the equality 
	$$
	(\alpha \cdot \theta) \cdot o = \theta \cdot (o \cdot \alpha), 
	$$
	one can obtain the relations (A7) and (A8).
	The first of this relations implies P3).
	
	\hspace{10mm}
	By the equality 
	$$
	\theta \cdot (o \cdot a) = \varepsilon \cdot a, \quad a \in M
	$$
the property P1(ii) and the relation (A9) are obtained. 
As consequence of (A9) we get that $\theta^{-1} = \Xi (o, o)$.
Finally, the property P1(i) follows immediately from the following Lemma.
		
		\hspace{10mm}
		{\bf Lemma 3.1.}  
		{\it For elements $x, a, b, \in M$ the relations 
		$[x, a] = b$ and $x \in M \cap H(b \cdot a^{-1})$ 
		are equivalent.} 
		
		(Note that $M \cap H(b \cdot a^{-1})$ has a unique element, since $M$ is a right  transversal to $H$.)
			
		\hspace{10mm}
		The Lemma 3.1 is true, because
		$$
		[x, a] = b 	 \, \Leftrightarrow \, x \cdot a = (x, a) \cdot [x, a] = (x, a) \cdot b  \, \Leftrightarrow \, 
		x  = (x, a) \cdot b \cdot a^{-1} \, \Leftrightarrow \, x \in M \cap H(b \cdot a^{-1}).
		$$

	\hspace{10mm}	
	Let $g: (G, H, M) \rightarrow (G', H', M')$ be a morphism of group triples, 
	$M_H$ and $M'_{H'}$ be the associated with the given triples hypergroups over the  group. 
	Then by restrictions of the corresponding group homomorphism $g: G \rightarrow G'$ 
	we get a group homomorphism $f_0: H \rightarrow H'$ 
	and a map of sets $f_1: M \rightarrow M'$ such that 
	$f = (f_0, f_1): M_H \rightarrow M'_{H'}$ is a morphism of hypergroups.

\hspace{10mm}	
Let ${\mathcal G}Trip_O$, ${\mathcal H}g_O$ be the classes of objects and 
${\mathcal G}Trip_M$, ${\mathcal H}g_M$ be the classes of morphisms 
of corresponding categories.
Consider the mappings
${\bf H}_O: {\mathcal G}Trip_O \rightarrow {\mathcal H}g_O$, 
which sends any group triple to the associated with him hypergroup over the group, and 
${\bf H}_M: {\mathcal G}Trip_M \rightarrow {\mathcal H}g_M$, 
which sends any morphism of triples to the corresponding morphism of associated hypergroups.
The pair ${\bf H} = ({\bf H}_O, {\bf H}_M)$ determines a (covariant) functor 
	from ${\mathcal G}Trip$ to ${\mathcal H}g$.


	\vspace{3mm}	
	
	{\bf 4. The exact product, associated with hypergroups over the group.}
	In this section the functor ${\bf T}: {\mathcal H}g \rightarrow {\mathcal G}Trip$ is defined.	
This definition is based on the construction of  exact product, associated with a hypergroup over the group.

\hspace{10mm}
Let $M_H$ be an arbitrary hypergroup over the group. 
Consider the set of all	two-letter words $\alpha a$, where $\alpha \in H$ and $a \in M$. 
Define the product of two such words by formula
$$
\alpha a \cdot \beta b = (\alpha \cdot {}^a\beta \cdot (a^\beta, b))[a^\beta, b].
$$
\hspace{10mm}
{\bf Proposition 4.1.} {\it The set of all two-letter words together with the above defined binary operation 
forms a group.}

\hspace{10mm}
This group is called {\it the exact product, associated wit the hypergroup over the group $M_H$}, 
and is denoted $H \odot M$.

\hspace{10mm}
{\bf Proof.} 
The associative law of the considered operation: 
$$
		(\alpha a \cdot \beta b) \cdot \gamma c = \alpha a \cdot (\beta b \cdot \gamma c)
$$
for any two-letter words $\alpha a, \beta b, \gamma c$, 
is true, because
$$
(\alpha a \cdot \beta b) \cdot \gamma c =
((\alpha \cdot {^a}\beta \cdot (a^\beta, b))\cdot {^{[a^\beta, b]}}\gamma \cdot ([a^\beta, b]^\gamma, \, c))
[[a^\beta, b]^\gamma, \, c]
$$
and
$$
\alpha a \cdot (\beta b \cdot \gamma c) =
(\alpha \cdot {^a}(\beta \cdot {^b}\gamma \cdot (b^\gamma, c)) \cdot
(a^{\beta \cdot {^b}\gamma \cdot (b^\gamma, c)}, \, [b^\gamma, c])
[a^{\beta \cdot {^b}\gamma \cdot (b^\gamma, c)}, \, [b^\gamma, c]].
$$
Here
$$
(\alpha \cdot {^a}(\beta \cdot {^b}\gamma \cdot (b^\gamma, c)) \cdot (a^{\beta \cdot {^b}\gamma \cdot (b^\gamma, c)}, \, [b^\gamma, c]) =
$$
$$
\alpha \cdot {^a}\beta \cdot {^{a^\beta}}({^b}\gamma) \cdot {^{a^{\beta \cdot {^b}\gamma}}}(b^\gamma, c) \cdot
(a^{\beta \cdot {^b}\gamma \cdot (b^\gamma, c)}, \, [b^\gamma, c]) =
$$
$$
\alpha \cdot {^a}\beta \cdot (a^\beta, \, b) \cdot {^{[a^\beta, \, b]}}\gamma
\cdot ((a^\beta)^{{^b}\gamma}, \, b^\gamma)^{-1}
 \cdot {^{a^{\beta \cdot {^b}\gamma}}}(b^\gamma, c) \cdot
(a^{\beta \cdot {^b}\gamma \cdot (b^\gamma, c)}, \, [b^\gamma, c]) =
$$
$$
\alpha \cdot {^a}\beta \cdot (a^\beta, \, b) \cdot {^{[a^\beta, \, b]}}\gamma
\cdot ([a^{\beta \cdot {^b}\gamma}, \, b^\gamma], \, c) =
$$
$$
((\alpha \cdot {^a}\beta \cdot (a^\beta, b))\cdot {^{[a^\beta, b]}}\gamma \cdot ([a^\beta, b]^\gamma, \, c))
$$
according to (A1), (A3), (A5), (A2), and
$$
[[a^\beta, b]^\gamma, \, c] =
[a^{\beta \cdot {^b}\gamma \cdot (b^\gamma, c)}, \, [b^\gamma, c]]
$$
according to (A2), (A4).

\hspace{10mm}
Similarly is checkes that $\theta o$ is a left neutral element for the considered operation.

\hspace{10mm}
Now we check  that for any two-letter words  $\alpha a, \beta b$ there exists a unique  
two-letter word $\xi x$ such that $\xi x \cdot \alpha a = \beta b$.
This equation is equivalent to the relation 
$$
(\xi \cdot {}^x\alpha \cdot (x^\alpha, a))[x^\alpha, a] = \beta b,
$$
and, consequently, to the pair of relations
$$
\xi \cdot {}^x\alpha \cdot (x^\alpha, a) = \beta, \quad [x^\alpha, a] = b.
$$
By using the denotation $b/a$ for the unique solution of the equation $[y, a ] = b$ 
we get that the last system of two equations has exactly one solution
$$
x = (b/a)^{\alpha^{-1}}, \quad \xi = \beta \cdot ({}^x\alpha \cdot (b/a, a))^{-1}.
$$ 
The Proposition 4.1 is proved.


		\hspace{10mm}
		{\bf Example 4.1.}
Let $M_H$ be a hypergroup over the group with trivial
structural mappings $\Phi, \Psi, \Lambda$, i.e. 
$$
a^\alpha = a, \quad {^a}\alpha = \alpha, \quad  (a, b) = \varepsilon
$$
for any $\alpha \in H, \, a, b \in M$.
Then $(M, \, \Xi)$ is a group and the  exact product, 
associated with the hypergroup $M_H$, 
coincides with the direct product of groups $H$ and $(M, \, \Xi)$.

		\hspace{10mm}
		{\bf Example 4.2.}
Suppose that only the structural mappings $\Phi$ and $\Lambda$ of the hypergroup $M_H$  are trivial.
Then again $(M, \, \Xi)$ is a group and the  exact product, associated with the
hypergroup $M_H$, is a semidirect product of $H$ by $(M, \, \Xi)$.

		\hspace{10mm}
		{\bf Example 4.3.}
If only the structural mapping  $\Lambda$ is trivial,
then  $(M, \, \Xi)$ is a group as well, and the corresponding exact product
is the {\it general product} of $H$ and $(M,\, \Xi)$ in the sense of B. Neumann
(\cite{N}, see also \cite{K}, p. 485).

\hspace{10mm}
Thus, the notion of the exact product, associated with a hypergroup over the group, 
generalizes the notions of the  direct product, of the semidirect product and of the  general product of two groups.


\hspace{10mm}
Let $M_H$ be an arbitrary hypergroup over the group and
${\overline G} = H \odot M$.
Consider the maps
$$
f_0: H \rightarrow {\overline G}, \quad \alpha \mapsto {\overline \alpha} = (\alpha \cdot \theta) o,
\quad \quad
f_1: M \rightarrow {\overline G}, \quad a \mapsto {\overline a} = \varepsilon a,
$$
and their images  ${\overline H} = {\rm im} \, f_0$, ${\overline M} = {\rm im} \, f_1$.

\hspace{10mm}
Then $\overline H$ is a subgroup of $\overline G$ isomorphic to $H$,
and 
the subset $\overline M$ of $\overline G$ is bijective to $M$ and  is a complementary set
to the subgroup $\overline H$.

\hspace{10mm}
Thus, there exists a canonical mapping for classes of objects of the corresponding categories
$$
{\bf T}_O: {\mathcal H}g_O \rightarrow {\mathcal G}Trip_O, \quad M_H \mapsto ({\overline G}, {\overline H}, {\overline M}).
$$

\hspace{10mm}
Now let $f = (f_0, f_1): M_H \rightarrow M'_{H'}$ be a morphism of hypergroups, 
$({\overline G}, {\overline H}, {\overline M})$ and  $({\overline {G'}}, {\overline {H'}}, {\overline {M'}})$ 
be the corresponding group triples.
Let the map ${\overline g}: {\overline G} \rightarrow {\overline G}'$ is defined by formula
$$
{\overline g} (\alpha a) = f_0(\alpha) f_1(a).
$$
 Then $\overline g$ is a group homomorphism and determines  
a morphism of group triples   
$$
{\overline g}: ({\overline G}, {\overline H}, {\overline M}) \rightarrow ({\overline {G'}}, {\overline {H'}}, {\overline {M'}}).
$$
Thus, there is a canonical mapping
$$
{\bf T}_M: {\mathcal H}g_M \rightarrow {\mathcal G}Trip_M, \quad 
f \mapsto {\overline g}.
$$

\hspace{10mm}
 The pair of mappings ${\bf T}  = ({\bf T} _O, {\bf T} _M)$ gives  
a (covariant) functor from ${\mathcal H}g$ to ${\mathcal G}Trip$.


	\vspace{3mm}	
	{\bf 5. Proof of   equivalence of categories ${\mathcal H}g$ and ${\mathcal G}Trip$.}
				The proof of the relation ${\bf T} \circ {\bf H}   \approx   1_{{\mathcal H}g}$ 
		is based on the following result.
		
		\hspace{10mm}
		{\bf Proposition 5.1.} 
		(The  universal property of the standard construction of hypergroups over the group.) 
		{\it  By the standard construction of hypergroups over the group, 
any hypergroup over the group (up to isomorphism) can be obtained from a group triple. 
More exactly, let $M_H$ be an arbitrary hypergroup over the group,
 $({\overline  G}, {\overline H}, {\overline M})$ be the group triple,  canonically associated with $M_H$,
${\overline M}_{\overline H}$ be the hypergroup over the group, obtained by the standard construction from  
the group triple  $({\overline  G}, {\overline H}, {\overline M})$. 
Then there exists a canonical isomorphism ${\overline f}: M_H \rightarrow {\overline M}_{\overline H}$.
}

	\hspace{10mm}	
	{\bf Proof.}
	Let ${\overline f} = ({\overline f}_0,  {\overline f}_1)$, where
	$$
{\overline f}_0: H \rightarrow {\overline H}, \quad \alpha \mapsto {\overline \alpha} = (\alpha \cdot \theta) o,
\quad \quad
{\overline f}_1: M \rightarrow {\overline M}, \quad a \mapsto {\overline a} = \varepsilon a
$$
are, respectively, an isomorphism of groups and a bijection of sets. 
To check the conditions $(M\Phi) - (M\Lambda)$ for the pair $({\overline f}_0, {\overline f}_1)$, 
we note that they have the following form:
$$
(M\Phi) \,\, {\overline a}^{\overline \alpha} = \overline{a^\alpha}, \quad 
(M\Psi) \,\, {}^{\overline a}{\overline \alpha} = \overline{{}^a\alpha}, \quad
(M\Xi)  \,\, [{\overline a}, {\overline b}] = {\overline {[a, b]}},  \quad
(M\Lambda) \,\,  ({\overline a}, {\overline b}) = {\overline {(a, b)}}.
$$

\hspace{10mm}
The lemma 5.1 is proved by a direct calculation.

\hspace{10mm}
{\bf Lemma 5.1.}
{\it For any elements $\xi \in H, \, x \in M$ ${\overline \xi} \cdot {\overline x} =  \xi x$.}

\hspace{10mm}	
According this Lemma  
$$
\overline{{}^a\alpha} \cdot \overline{a^\alpha} = {}^a\alpha a^\alpha, \quad 
\overline{(a, b)} \cdot \overline{[a, b]} = (a, b)[a, b],
$$
and using the definition of the multiplication in $\overline G$ and the relations (A1), (A11), (A10), P2(i), (A6)
we obtain
$$
{}^{\overline a}{\overline \alpha} \cdot{\overline a}^{\overline \alpha} = {\overline a} \cdot{\overline \alpha} =
(\varepsilon a) \cdot (\alpha \cdot \theta)o = 
(\varepsilon \cdot {}^a(\alpha \cdot \theta) \cdot (a^{\alpha \cdot \theta}, o))[a^{\alpha \cdot \theta}, o] = 
({}^a\alpha \cdot {}^{a^\alpha}\theta \cdot {}^{a^{\alpha \cdot \theta}}(\theta^{-1}))((a^{\alpha \cdot \theta})^{\theta^{-1}}) =  
$$
$$
({}^a\alpha \cdot {}^{a^\alpha}(\theta \cdot \theta^{-1}))a^\alpha  = {}^a\alpha a^\alpha . 
$$
	and similarly 
	$$
	(\overline{a}, \overline{b}) \cdot [\overline{a}, \overline{b}] = \overline{a} \cdot  \overline{b} = \varepsilon a \cdot \varepsilon b = 
	(\varepsilon \cdot {}^a\varepsilon \cdot (a^\varepsilon, b))[a^\varepsilon, b] = (a, b)[a, b].
	$$

\hspace{10mm}
{\bf Corollary 5.1.1.} 
{\it By the system of isomorphisms 
$$
{\overline f}: M_H \rightarrow {\overline M}_{\overline H}, \quad  M_H \in {\mathcal H}g_O
$$
an isomorphism of functors $ 1_{{\mathcal H}g}  \approx {\bf T} \circ {\bf H}$ is determined.} 

\hspace{10mm}
{\bf Proof.}
By a direct calculation it is proved that 
for any morphism of hypergroups $f: M_N \rightarrow M'_{H'}$ we have a commutative square 
\begin{center}
$M_H \longrightarrow M'_{H'}$ 

${\overline f} \downarrow \hspace{12mm}  \downarrow {\overline {f'}} \hspace{2mm}$

${\overline M}_{\overline H} \longrightarrow {\overline {M'}}_{\overline {H'}}$.
\end{center}
(Here, the upper  and  lower horizontal arrows represent the morphisms $f$ and $({\bf T} \circ {\bf H})(f)$, respectively.)

\hspace{10mm}
Similarly  a natural isomorphism $1_{{\mathcal G}Trip} \rightarrow {\bf H} \circ {\bf T}$ is constructed..


\hspace{10mm}
{\bf Proposition 5.2.} (The universal property of the exact product, associated with hypergroups over the group.)
{\it Any group triple is isomorphic to a group triple, associated with a hypergroup over the group.
More exactly, let $(G, H, M)$ be an arbitrary group triple, $M_H$ be the hypergroup, obtauned by the standard construction from $(G, H, M)$, 
$(\overline{G}, \overline{H}, \overline{M})$ be the exact product, associated with $M_H$. 
Let the map $\overline{g}: G \rightarrow {\overline G}$ is defined by $\overline{g}(x) = \alpha a$, where $x \in G, x = \alpha \cdot a, \alpha \in H, a \in M$. 
This map $\overline g$ is a group homomorphism and determines a canonical isomorphism of group triples 
$(G, H, M) \rightarrow (\overline{G}, \overline{H}, \overline{M})$.}

\hspace{10mm}
{\bf Corollary 5.2.1.}
{\it Let $g: (G, H, M) \rightarrow (G', H', M')$ be an arbitrary morphism of group triples, 
while $\overline{g}: (G, H, M) \rightarrow ({\overline G}, {\overline H}, {\overline M})$ and 
$\overline{g'}: (G', H', M') \rightarrow (\overline{G'}, \overline{H'}, \overline{M'})$ be 
the canonical isomorphisms, associated with triples $(G, H, M)$, $(G', H', M)$ and 
${\tilde g} = ({\bf H} \circ {\bf T})(g)$. 
Then $g \circ {\overline{g'}} = {\overline g} \circ {\tilde g}$.} 

This terminate the proof of the equivalence of categories ${\mathcal H}g$ and ${\mathcal G}Trip$.


\end{document}